
%
\newlinechar`\^^J%
\message{^^J Attention, vous compilez un fichier comportant de nombreuses ^^J20
macros specifiques a l'Institut Fourier, il est recommande de lire ^^J
les commentaires figurant au debut de ce fichier. Ce sera  ^^J
necessaire si la compilation vous indique des erreurs.^^J^^J}
\message{ Please learn french and read the included comments in case of trouble.^^J
This file should compile with *plain* TeX. It inputs two files from the ^^J
AMS TeX distribution and requires AMS fonts: euler, symbol and extracm.^^J^^J}
%
%
\input amssym.def
\input amssym 
\catcode`\|=13 
\font\helvb=cmssbx10
\font\eightrm=cmr8
\font\eighti=cmmi8
\font\eightsy=cmsy8
\font\eightbf=cmbx8
\font\eighttt=cmtt8
\font\eightit=cmti8
\font\eightsl=cmsl8
\font\sixrm=cmr6
\font\sixi=cmmi6
\font\sixsy=cmsy6
\font\sixbf=cmbx6
\skewchar\eighti='177 \skewchar\sixi='177
\skewchar\eightsy='60 \skewchar\sixsy='60

\def\tenpoint{%
  \textfont0=\tenrm \scriptfont0=\sevenrm \scriptscriptfont0=\fiverm
  \def\rm{\fam0\tenrm}%
  \textfont1=\teni \scriptfont1=\seveni \scriptscriptfont1=\fivei
  \def\mit{\fam\@ne}\def\oldstyle{\fam1\teni}%
  \textfont2=\tensy \scriptfont2=\sevensy \scriptscriptfont2=\fivesy
    \def\itfam{4}\textfont\itfam=\tenit
  \def\it{\fam\itfam\tenit}%
  \def\slfam{5}\textfont\slfam=\tensl
  \def\sl{\fam\slfam\tensl}%
  \def\bffam{6}\textfont\bffam=\tenbf \scriptfont\bffam=\sevenbf
  \scriptscriptfont\bffam=\fivebf
  \def\bf{\fam\bffam\tenbf}%
  \def\ttfam{7}\textfont\ttfam=\tentt
  \def\tt{\fam\ttfam\tentt}%
  \abovedisplayskip=6pt plus 2pt minus 6pt
  \abovedisplayshortskip=0pt plus 3pt
  \belowdisplayskip=6pt plus 2pt minus 6pt
  \belowdisplayshortskip=7pt plus 3pt minus 4pt
  \smallskipamount=3pt plus 1pt minus 1pt
  \medskipamount=6pt plus 2pt minus 2pt
  \bigskipamount=12pt plus 4pt minus 4pt
  \normalbaselineskip=12pt
  \setbox\strutbox=\hbox{\vrule height8.5pt depth3.5pt width0pt}%
  \normalbaselines\rm}

\def\eightpoint{%
  \textfont0=\eightrm \scriptfont0=\sixrm \scriptscriptfont0=\fiverm
  \def\rm{\fam0\eightrm}%
  \textfont1=\eighti \scriptfont1=\sixi \scriptscriptfont1=\fivei
  \def\oldstyle{\fam1\eighti}%
  \textfont2=\eightsy \scriptfont2=\sixsy \scriptscriptfont2=\fivesy
  \textfont\slfam=\eightit
  \def\sl{\fam\itfam\eightit}%
  \textfont\slfam=\eightsl
  \def\sl{\fam\slfam\eightsl}%
  \textfont\bffam=\eightbf \scriptfont\bffam=\sixbf
  \scriptscriptfont\bffam=\fivebf
  \def\bf{\fam\bffam\eightbf}%
  \textfont\ttfam=\eighttt
  \def\tt{\fam\ttfam\eighttt}%
  \abovedisplayskip=9pt plus 2pt minus 6pt
  \abovedisplayshortskip=0pt plus 2pt
  \belowdisplayskip=9pt plus 2pt minus 6pt
  \belowdisplayshortskip=5pt plus 2pt minus 3pt
  \smallskipamount=2pt plus 1pt minus 1pt
  \medskipamount=4pt plus 2pt minus 1pt
  \bigskipamount=9pt plus 3pt minus 3pt
  \normalbaselineskip=9pt
  \setbox\strutbox=\hbox{\vrule height7pt depth2pt width0pt}%
  \normalbaselines\rm}

\font\petcap=cmcsc10

\tenpoint
\hsize=12.5cm
\vsize=19cm
\parskip 5pt plus 1pt
\parindent=1cm
\baselineskip=13pt
\hoffset=-0.1cm 
\def\footnoterule{\kern-6pt
  \hrule width 2truein \kern 5.6pt} 

\def\ie{{\sl i.e.\ }}

\def\omini{\raise 1ex\hbox{\ept o}}
\def\emini{\raise 1ex\hbox{\ept e}}
\def\ermini{\raise 1ex\hbox{\ept er}}
\def\remini{\raise 1ex\hbox{\ept re}}

\def\lead{\leaders\hbox to 10pt{\hss.\hss}\hfill}
\def\somt#1|#2|{\vskip 8pt plus1pt minus 1pt
                \line{#1\lead #2}}
\def\soms#1|#2|{\vskip 2pt
                \line{\qquad #1\lead #2}} 
\def\somss#1|#2|{\vskip 1pt  
                 \line{\qquad\qquad #1\lead #2}}


\def\aujour{\ifnum\day=1 1\ermini\else\number\day\fi\
\ifcase\month\or janvier\or f\'evrier\or mars\or avril\or mai\or juin\or
juillet\or aout\or septembre\or octobre\or novembre\or d\'ecembre\fi\
\number\year}
\def\today{\ifcase\month\or january \or february \or march \or april
\or may \or june\or july\or august \or september\or october\or november\or
december\fi\ \number\day , \number\year}

\newskip\afterskip
\catcode`\@=11
\def\p@int{.\par\vskip\afterskip\penalty100} 
\def\p@intir{\discretionary{.}{}{.\kern.35em---\kern.7em}}
\def\pointir{\afterassignment\pointir@\global\let\next=}
\def\pointir@{\ifx\next\par\p@int\else\p@intir\fi\egroup\next}
\catcode`\@=12
\def|{\relax\ifmmode\vert\else\findef\fi}
\def\findef{\errhelp{Cette barre verticale ne correspond ni a un \vert
mathematique
                        ni a une fin de definition, le contexte doit vous
indiquer ce qui manque.
                        Si vous vouliez inserer un long tiret, le codage
recommande est ---,
                        dans tous les cas, la barre fautive a ete supprimee.}%
                        \errmessage{Une barre verticale a ete trouvee en
mode texte}}

\def\TITR#1|{\null{\mss\baselineskip=17pt
                           \vskip 3.25ex plus 1ex minus .2ex
                           \leftskip=0pt plus \hsize
                           \rightskip=\leftskip
                           \parfillskip=0pt
                           \noindent #1
                           \par\vskip 2.3ex plus .2ex}}

\def\auteur#1|{\penalty 500
               \vbox{\centerline{\si
                 \iffrance par
\fi #1}
                \vskip 10pt}\penalty 500}

\def\resume#1|{\penalty 100
                           {\leftskip=\parindent
                            \rightskip=\leftskip
                            \eightpoint\bgroup\petcap \skip\afterskip=0pt
                             \iffrance R\'esum\'e \else Abstract \fi\pointir
                            #1 \par}
                           \penalty -100}

\def\titre#1|{\null\baselineskip14pt
                           {\helvb
                           \vskip 3.25ex plus 1ex minus .2ex
                           \leftskip=0pt plus \hsize
                           \rightskip=\leftskip
                           \parfillskip=0pt
                           \noindent #1
                           \par\vskip 2.3ex plus .2ex}}


\def\section#1|{
                                \bgroup\bf
                                 \par\penalty -500
                                 \vskip 3.25ex plus 1ex minus .2ex
                                 \skip\afterskip=1.5ex plus .2ex
                                  #1\pointir}

\def\ssection#1|{
                                 \bgroup\petcap
                                  \par\penalty -200
                                  \vskip 3.25ex plus 1ex minus .2ex
                          \skip\afterskip=1.5ex plus .2ex
                                   #1\pointir}


\def\th#1|{
                   \bgroup \sl
                        \def\findef{\egroup\par}
                        \bgroup\petcap
                         \par\vskip 2ex plus 1ex minus .2ex
                         \skip\afterskip=0pt
                           #1\pointir}

\def\defi#1|{
                   \bgroup \rm
                        \def\findef{\egroup\par}
                        \bgroup\petcap
                         \par\vskip 2ex plus 1ex minus .2ex
                         \skip\afterskip=0pt
                           #1\pointir}

\def\rque#1|{\bgroup \sl
                          \par\vskip 2ex plus 1ex minus .2ex\skip\afterskip=0pt
                          #1\pointir}

\def\dem{\bgroup \sl
                  \par\vskip 2ex plus 1ex minus .2ex\skip\afterskip=0pt
                  \iffrance D\'emonstration\else Proof\fi\pointir}

\def\preuve{\bgroup \sl
                  \par\vskip 2ex plus 1ex minus .2ex\skip\afterskip=0pt
                  \iffrance Preuve\else Proof\fi\pointir}


\def\_#1{_{\baselineskip=.7 \baselineskip
                                                                      
\vtop{\halign{\hfil$\scriptstyle{##}$\hfil\cr #1\crcr}}}}

\def\build#1#2\fin{\mathrel{\mathop{\kern0pt#1}\limits#2}}

\def\frac#1/#2{\leavevmode\kern.1em
   \raise.5ex\hbox{$\scriptstyle #1$}\kern-.1em
      /\kern-.15em\lower.25ex\hbox{$\scriptstyle #2$}}

{\obeylines
}

\newif\ifchiffre
\def\chiffre{\chiffretrue}
\chiffre
\newdimen\laenge
\def\lettre#1|{\setbox3=\hbox{#1}\laenge=\wd3\advance\laenge by 3mm
\chiffrefalse}
\def\article#1|#2|#3|#4|#5|#6|#7|%
    {{\ifchiffre\leftskip=7mm\noindent
     \hangindent=2mm\hangafter=1
\llap{[#1]\hskip1.35em}\bgroup\petcap #2\pointir {\sl #3}, {\rm #4}
\nobreak{\bf #5}
({\oldstyle #6}), \nobreak #7.\par\else\noindent \advance\laenge by 4mm
\hangindent=\laenge\advance\laenge by -4mm\hangafter=1
\rlap{[#1]}\hskip\laenge\bgroup\petcap #2\pointir {\sl #3}, #4 {\bf #5}
({\oldstyle
#6}), #7.\par\fi}} 
\def\livre#1|#2|#3|#4|#5|%
    {{\ifchiffre\leftskip=7mm\noindent
    \hangindent=2mm\hangafter=1
\llap{[#1]\hskip1.35em}\bgroup\petcap #2\pointir{\sl #3}, #4, {\oldstyle
#5}.\par
\else\noindent
\advance\laenge by 4mm \hangindent=\laenge\advance\laenge by -4mm
\hangafter=1
\rlap{[#1]}\hskip\laenge\bgroup\petcap #2\pointir
{\sl  #3}, #4, {\oldstyle #5}.\par\fi}}
\def\divers#1|#2|#3|#4|%
    {{\ifchiffre\leftskip=7mm\noindent
    \hangindent=2mm\hangafter=1
     \llap{[#1]\hskip1.35em}\bgroup\petcap #2\pointir #3, {\oldstyle #4}.\par
\else\noindent
\advance\laenge by 4mm \hangindent=\laenge\advance\laenge by -4mm
\hangafter=1
\rlap{[#1]}\hskip\laenge\bgroup\petcap #2\pointir #3,{\oldstyle #4}.\par\fi}}
\def\div#1|#2|#3|#4|
{{\ifchiffre\leftskip=7mm\noindent
\hangindent=2mm\hangafter=1
\llap{[#1]\hskip1.35em}\bgroup\petcap #2\pointir {\sl  #3},{\oldstyle #4}.\par
\else\noindent
\advance\laenge by 4mm \hangindent=\laenge\advance\laenge by -4mm
\hangafter=1
\rlap{[#1]}\hskip\laenge\bgroup\petcap #2\pointir {\sl  #3},{\oldstyle
#4}.\par\fi}}


\font\si=cmssi10

\font\mss=cmss12 scaled \magstep1 

\hsize=12.5cm
\vsize=19cm
\parskip 5pt plus 1pt
\parindent=1cm
\baselineskip=13pt

  \hoffset=0.7cm
 \voffset=0.8cm

\font\cms=cmbsy10 at 5pt
\font\gros=cmex10 at 17.28pt

\newif\iffrance

\def\TITR#1|{\null{\mss\baselineskip=17pt
                           \vskip 3.25ex plus 1ex minus .2ex
                           \leftskip=0pt plus \hsize
                           \rightskip=\leftskip
                           \parfillskip=0pt
                           \noindent #1
                           \par\vskip 2.3ex plus .2ex}}

\def\auteur#1|{\penalty 500
               \vbox{\centerline{\si
                 \iffrance par 
\fi #1}
                \vskip 10pt}\penalty 500}

\def\resume#1|{\penalty 100
                           {\leftskip=\parindent
                            \rightskip=\leftskip
                            \eightpoint\bgroup\petcap \skip\afterskip=0pt
                             \iffrance R\'esum\'e \else Abstract \fi\pointir
                            #1 \par}
                           \penalty -100}

\def\titre#1|{\null\penalty-500\baselineskip14pt
                           {\helvb
                           \vskip 3.25ex plus 1ex minus .2ex
                           \leftskip=0pt plus \hsize
                           \rightskip=\leftskip
                           \parfillskip=0pt
                           \noindent #1
                           \par\nobreak\vskip 2.3ex plus .2ex}\penalty 5000}


\def\prodd{\mathop{\prod}\limits}
\def\summ{ \mathop{\sum}\limits}

\def\opp{\mathop{\bigoplus}\limits}

\def\hfl{{\hbox to 12mm{\rightarrowfill}}}

\def\egal{{{\lower5pt\hbox to 0.5cm{$\hrulefill$}}\atop
{\raise10pt\hbox to 0.5cm{$\hrulefill$}}}}

\def\upuparrow{\eqalign{&\big\uparrow\cr\noalign{\vskip-14pt}&\kern-1.8pt
               \uparrow\cr}}
\def\ddarrow{\eqalign{&\raise15pt\hbox{$\big\downarrow$}\cr
             \noalign{\vskip-25pt} &\big\downarrow\cr}}

\def\cu{\hbox{\cms\char'133}}
\def\ca{\hbox{\cms\char'134}} 
 \def\hookup{{\lower6.3pt\hbox{\cu}\kern-4.32pt\big\uparrow}}
 \def\hookdown{\raise5.3pt\hbox{\ca}\kern-4.32pt\lower3pt\hbox{\big\downarrow}}

\def\grtilde{\kern10pt\hbox{\gros\char'147}}

\let\le=\leq 
\let\ge=\geq 

\def\S{\mathhexbox278\kern.15em}
\def\\ {\smallsetminus}

\def\wt{\widetilde}

\def\ept{\eightpoint}

\def\la{\longrightarrow}
\def\ld {, \ldots,}

\def\coin{\mathrel{\raise2pt\hbox{$\scriptstyle
|$}\kern-2pt\lower2pt\hbox{$-$}}}
\def\rect{\mathrel{\lower2pt\hbox{$-$}\kern-2pt\raise2pt\hbox{$\scriptstyle
|$}}}


{\obeylines
}


\long\def\nomm#1|#2|#3|{\line{$\vtop{\hsize2.5cm #1~}\vtop{\hsize1cm #2~}
\vtop{\hsize=9cm\normalbaselines\parshape 1 0cm 9cm #3.}$}
\medskip}

\long\def\nom#1|#2|{\centerline{$\vtop{\hsize=3cm{\bf #1}~: }
\vtop{\hsize=9.5cm\normalbaselines\parshape 1 0cm 9.5cm #2.}$}
\medskip}

\def\boxit#1#2{\hbox{\vrule
 \vbox{\hrule\kern#1
  \vtop{\hbox{\kern#1 #2\kern#1}%
   \kern#1\hrule}}%
 \vrule}}

\newbox\texte
\def\texteencadre#1|#2|#3|{\setbox1=\vbox{#3}
\setbox\texte=\vbox{\hrule height#1pt%
\hbox{\vrule width#1pt\kern#2pt\vbox{\kern#2pt \hbox{\box1}\kern#2pt}%
\kern#2pt\vrule width#1pt}\hrule height#1pt}}

\def\blanc{\hbox to 0pt{\vrule height13pt depth5pt width0pt}}
\def\blancs{\hbox to 0pt{\vrule height9pt depth5pt width0pt}}
\def\blan{\hbox to 0pt{\vrule height6pt depth5pt width0pt}}
\def\bla{\hbox to 0pt{\vrule height7pt depth5pt width0pt}}
\def\hrulefill{\leaders\hrule height0.2pt\hfill}

\def\bb#1&#2\cr{
\hbox to 1cm{\hfil\strut#1~\vrule}
\hbox to 1cm{\hfil\strut#2~\vrule}
}
\def\aa#1&#2&#3&#4&#5&#6&#7&#8&#9\cr{\line{$
\hbox to 3cm{\vrule width 1pt \strut~#1\hfil\vrule}
\hbox to 1cm{\hfil\strut#2~\vrule}
\hbox to 1cm{\hfil\strut#3~\vrule}
\hbox to 1cm{\hfil\strut#4~\vrule}
\hbox to 1cm{\hfil\strut#5~\vrule}
\hbox to 1cm{\hfil\strut#6~\vrule}
\hbox to 1cm{\hfil\strut#7~\vrule}
\hbox to 1cm{\hfil\strut#8~\vrule width 1pt}
\bb#9\cr
$\hfil}}


           \def\cb{{\Bbb C}}

           \def\hb{{\Bbb H}}

\def\oc{{\cal O}}                      
                      
           \def\qb{{\Bbb Q}}            
           \def\rb{{\Bbb R}}

           \def\zb{{\Bbb Z}}            


                      \def\ggg{{\goth g}}

\font\tenmib=cmmib10 
\font\sevenmib=cmmib7
\font\fivemib=cmmib5
\expandafter\chardef\csname pre boldmath.tex at\endcsname=\the\catcode`\@
\catcode`\@=11
\def\hexanumber@#1{\ifcase#1 0\or 1\or 2\or 3\or 4\or 5\or 6\or 7\or 8\or
 9\or A\or B\or C\or D\or E\or F\fi}

\skewchar\tenmib='177 \skewchar\sevenmib='177 \skewchar\fivemib='177

\newfam\mibfam   
\textfont\mibfam=\tenmib \scriptfont\mibfam=\sevenmib
\scriptscriptfont\mibfam=\fivemib
\def\mib@hex{\hexanumber@\mibfam}
\mathchardef\bfalpha="0\mib@hex 0B
\mathchardef\bfbeta="0\mib@hex 0C
\mathchardef\bfgamma="0\mib@hex 0D
\mathchardef\bfdelta="0\mib@hex 0E
\mathchardef\bfepsilon="0\mib@hex 0F
\mathchardef\bfzeta="0\mib@hex 10
\mathchardef\bfeta="0\mib@hex 11
\mathchardef\bftheta="0\mib@hex 12
\mathchardef\bfiota="0\mib@hex 13
\mathchardef\bfkappa="0\mib@hex 14
\mathchardef\bflambda="0\mib@hex 15
\mathchardef\bfmu="0\mib@hex 16
\mathchardef\bfnu="0\mib@hex 17
\mathchardef\bfxi="0\mib@hex 18
\mathchardef\bfpi="0\mib@hex 19
\mathchardef\bfrho="0\mib@hex 1A
\mathchardef\bfsigma="0\mib@hex 1B
\mathchardef\bftau="0\mib@hex 1C
\mathchardef\bfupsilon="0\mib@hex 1D
\mathchardef\bfphi="0\mib@hex 1E
\mathchardef\bfchi="0\mib@hex 1F
\mathchardef\bfpsi="0\mib@hex 20
\mathchardef\bfomega="0\mib@hex 21
\mathchardef\bfvarepsilon="0\mib@hex 22
\mathchardef\bfvartheta="0\mib@hex 23
\mathchardef\bfvarpi="0\mib@hex 24
\mathchardef\bfvarrho="0\mib@hex 25
\mathchardef\bfvarsigma="0\mib@hex 26
\mathchardef\bfvarphi="0\mib@hex 27
\mathchardef\bfimath="0\mib@hex 7B
\mathchardef\bfjmath="0\mib@hex 7C
\mathchardef\bfell="0\mib@hex 60
\mathchardef\bfwp="0\mib@hex 7D
\mathchardef\bfpartial="0\mib@hex 40
\mathchardef\bfflat="0\mib@hex 5B
\mathchardef\bfnatural="0\mib@hex 5C
\mathchardef\bfsharp="0\mib@hex 5D
\mathchardef\bftriangleleft="2\mib@hex 2F
\mathchardef\bftriangleright="2\mib@hex 2E
\mathchardef\bfstar="2\mib@hex 3F
\mathchardef\bfsmile="3\mib@hex 5E
\mathchardef\bffrown="3\mib@hex 5F
\mathchardef\bfleftharpoonup="3\mib@hex 28
\mathchardef\bfleftharpoondown="3\mib@hex 29
\mathchardef\bfrightharpoonup="3\mib@hex 2A
\mathchardef\bfrightharpoondown="3\mib@hex 2B
\mathchardef\bflhook="3\mib@hex 2C 
\mathchardef\bfrhook="3\mib@hex 2D 
\mathchardef\bfldotp="6\mib@hex 3A 
\catcode`\@=\csname pre boldmath.tex at\endcsname

\overfullrule=0pt
\long\def\InsertFig#1 #2 #3 #4\EndFig{
\hbox{\hskip #1 mm$\vbox to #2 mm{\vfil\includegraphics{#3}}#4$}}
\long\def\LabelTeX#1 #2 #3\ELTX{\rlap{\kern#1mm\raise#2mm\hbox{#3}}}

\def\Spin{\mathop{\rm Spin}\nolimits}

\def\slg{{\goth s l}}

\def\cont{\mathop{\rm cont}\nolimits}
\def\Iso{\mathop{\rm Iso}\nolimits}
\def\Hom{\mathop{\rm Hom}\nolimits}

\def\oqb{\overline \qb}
\def\Ka{\hbox{Kazhdan}}

{\parindent=0pt\ept\footnote{\hphantom{$\!$}}{\parindent=0pt
{\sl Mots-cl{\'e}s}~: K{\"a}hler groups, property $T$.}}

\TITR STRUCTURE OF K{\"A}HLER GROUPS, I~: SECOND COHOMOLOGY|

\auteur Alexander REZNIKOV|
\centerline{April 1, 1998}

\section 0. Introduction|

Fundamental groups of complex projective varieties are very difficult
to understand. There is a tremendous gap between few computed examples
and few general theorems. The latter all deal with either linear finite
dimensional representations ([Sim]) or actions on trees ([Gr-Sch]); besides,
one knows almost nothing.

This paper presents a new general theorem, partially settling a
well-known conjecture of Carlson-Toledo (cf [Ko]).

\th Main Theorem|Let $\Gamma $ be a fundamental group of a compact
K{\"a}hler manifold. Assume $\Gamma $ is not \Ka. Then $H^2(\Gamma
,\rb) \ne 0$. Morover, let   $\Delta $ be a finitely presented group which is not  \Ka  and let     $\Gamma \to \Delta$ be a central extension. If  $\Gamma$ is a fundamental group of a  compact
K{\"a}hler manifold, then the natural map in the second real cohomology  $ H^2(\Delta ,\rb) \to   H^2(\Gamma ,\rb)  $ is nonzero.           |

\th Corollary|Suppose $\Delta $ is not \Ka~ .Let $\wt \Delta $ be  the universal central extension $0\to
H_2(\Delta ,\zb) \to \wt \Delta  \to \Delta  \to 1$. Then $\Gamma $ is
not a fundamental group of a compact
K{\"a}hler manifold.|

\rque Examples|

1. Any amenable group is not \Ka, so their universal central extensions are not K{\"a}hler.

2. Any lattice in $SU(n,1)$ is not \Ka. However, in this case one can do better, see remark in
section 8 below. It follows from the Corollary to the Main Theorem that a universal central extension of such a lattice is not a K{\"a}hler group. In contrast with this, any central extension with finite cyclic group as a center, whose extension class is 
a reduction $mod \ n$ of the K{\"a}hler class {\it is}  a  K{\"a}hler  group, as shown by Deligne, K\`ollar and Catanese [Ko].

\th Theorem 0.1|Let $\Gamma $ be a fundamental group of a complex
projective variety. Suppose $\Gamma $ has a Zariski dense rigid
representation in $SO(2,n)$, $n$ odd. Then

(i) $H^2(\Gamma ,\rb) \ne 0$.

(ii) Moreover, $H^2_b(\Gamma ,\rb) \ne 0$ and the canonical map $H^2_b
\to H^2$ is not zero.|

\th Corollary 0.1|Let $\Delta $ be a lattice in $SO(2,n)$, $n$ odd,
uniform or not.
Let $\wt \Delta $ be a universal central extension of $\Delta $. Then
$\Delta $ is not a fundamental group of a complex projective variety.|

\th Theorem 0.2|Let $\Gamma $ be a fundamental group of a complex
projective variety. Suppose $\Gamma $ has a Zariski dense rigid
representation in $Sp(4)$. Then

(i) $H^2(\Gamma ,\rb) \ne 0$.

(ii) Moreover, $H^2_b(\Gamma ,\rb) \ne 0$ and the canonical map $H^2_b
\to H^2$ is not zero.|

\th Corollary 0.2|Let $\Delta $ be a lattice in $Sp(4)$. Then $\wt
\Delta $ is not a fundamental group of a complex projective variety.|

\th Corollary 0.3|Lattices in $\Spin (2,n)$, $n$ odd, and
$\wt{Sp(4)}$ are not fundamental groups of a complex projective
variety.|

The next result shows that three-manifold groups wich are rich in the sence of [Re2] are not  K{\"a}hler .

\th Theorem 10.3|(Rich Three-manifold groups are not  K{\"a}hler ). Let $M^3$ be an irreducible atoroidal three-manofold. Suppose there exists a
Zariski dense homomorphism $\rho  : \pi _1(M) \to SL_2 (\cb)$. Then
$\Gamma  = \pi _1(M)$ is not K{\"a}hler.| 

\section 1. A geometric picture for rigid representations|

Let $Y$ be a compact K{\"a}hler manifold. All rigid irreducible
representations $\rho  : \pi _1(Y) \to SL(N,\cb)$ are conjugate to
representations landing in a $SU(m,n) \subset SL(N,\cb)$ with $m+n =
N$ ([Sim~1]) and have a structure of complex variation of Hodge
structure ([Sim~1]). Moreover, we can always arrange that this conjugate
representation is defined over $\oqb$ (see e.g. [Re~1]). Relabeling ,
we assume that $\rho $ itself is defined over $\oqb$. We assume that
moreover, $\rho $ is defined over $\oc(\oqb)$; by a conjecture
of Carlos Simpson ([Sim~1]) this is always the case. Let $\{\rho _i\}$ be all the
Galois twists of $\rho $, then $\rho _i$ are rigid therefore land in
$SU(m_i,n_i)$. The image of $\pi _1(Y)$ in $\prodd_i SU(m_i,n_i)$ is
discrete; call it $\Gamma $. 

Coming back to $\rho $, consider a corresponding $\theta $-bundle $E$
([Sim~1]). It has the following structure: $E = \opp_{p+q=k}
L^{p,q}$ and $\theta $ maps $L^{p,q}$ to $L^{p-1,q+1} \otimes \Omega
^1$. Any $\theta $-invariant subbundle of $E$ has negative degree; in
particular, the degree of $L^{k,0}$ is positive. A harmonic metric $K$
in $E$ is the unique metric satisfying the equation ([Hit]) $F_K =
[\theta ,\theta ^*]$. The hermitian connection $\nabla_K$ leaves all
$E^{p,q}$ invariant. The connection $\nabla_K + \theta +\theta ^*$ is
flat with monodromy $\rho $. Let $V$ be the corresponding flat
holomorphic bundle. In $V$, we have a flag of holomorphic subbundles
$F^p = V^{k,0} \oplus \cdots \oplus V^{p,k-p}$, where $V^{p,q}$ are
$E^{p,q}$ thought as $C^\infty$-subbundles of $V$ with its new
holomorphic structure. We have therefore a $\rho $-equivariant map $\wt Y
\build \la^s\fin D$, where $D$ is a corresponding Griffiths domain
([G], ch.~I--II). Changing the sign of $K$ alternatively on
$V^{p,q}$ we obtain a flat pseudo-hermitian metric in $V$.

So if $m = \opp_{p ~\hbox{\ept even}} \dim V^{p,q}$, $n = \opp_{p~
\hbox{\ept odd}} \dim V^{p,q}$, then $\rho $ lands in $SU(m,n)$. The
Griffiths domain  $D$ carries a horizontal distribution defined by
the condition that the derivative of $F_p$ lies in $F_{p+1}$. The
developing map $s$ is horizontal. Differentiating this condition we
obtain a second order equation ([Sim~1]) $[\theta ,\theta ] = 0$, in
other words for $Z,W \in T_xY$, $\theta (Z)$ and $\theta (W)$ commute.

Since the image of $\pi _1(Y)$ in $\prodd_i SU(m_i,n_i)$ is discrete,
we obtain

\th Proposition 1.1 {\rm (Geometric picture for rigid
representations)}|Let $\rho  : \pi _1(Y) \to SL(N,\cb)$ be a rigid
irreducible representation, defined over $\oc(\oqb)$. Then there exist
Griffiths domains $D_i = SU(m_i,n_i)/K_i$, a discrete group $\Gamma $
in $\prodd SU(m_i,n_i)$ and a horizontal holomorphic map
$$S : Y\la \prodd D_i/\Gamma $$
which induces $\rho $ and all its Galois twists.|

\rque Remark|Though a Griffiths domain $D$ is topologically a
fibration over a hermitian symmetric space with fiber a flag variety,
generally it does not have a $SU(m,n)$-invariant K{\"a}hler metric. So
the complex manifold $\prodd D_i/\Gamma $ is not K{\"a}hler.

\rque Remark|This proposition tells us that one cannot expect too many
compact K{\"a}hler manifolds to have a nontrivial linear representation of
their fundamental group, of finite dimension.

\th Lemma 1.2 {\rm (Superrigidity)}|

(1) Let $X = \Gamma \\ SU(m,n)/S(U(m) \times U(n)$ be  a compact
hermitian locally symmetric space of  Siegel type
$I$. Let $Y$ be  compact K{\"a}hler and let $f : Y\to X$ be continuous. If
$f_* : \pi _1(Y) \to SU(n,m)$ is rigid and Zariski dense, e.g. $f_* : \pi _1(Y) \to
\Gamma $ an isomorphism, then either $f$ is homotopic to a 
holomorphic map, or there exists a compact complex analytic space
$Y'$, $\dim Y' < \dim X$, and a holomorphic map $\varphi  : Y \to Y'$
such that $f$ is homotopic to a composition $Y \build \la^\varphi \fin
Y'\build \la^{f_1}\fin X$.

(2) Let $X = \Gamma \\ SO(2,n)/S(O(2)\times O(n))$ be of Siegel type
IV. Let $Y$
be compact K{\"a}hler and let $f : Y\to X$ be continuous. If $f_* : \pi
_1(Y) \to SO(2,n)$
is rigid and Zariski dense  then either $f$ is homotopic to a
holomorphic map $f_0$, or $n$ is even and  there exists a compact complex analytic space
$Y'$, $\dim Y' < \dim X$, and a holomorphic map $\varphi  : Y \to Y'$
such that $f$ is homotopic to a composition $Y \build \la^\varphi \fin
Y'\build \la^{f_1}\fin X$.

(3) Let $X = \Gamma \\ Sp(2n)/U(n)$ be of Siegel type III. Let $Y$ be
compact K{\"a}hler and let $f : Y \to X$ be continuous. If $f_* : \pi
_1(Y) \to  Sp(2n)$ is rigid and Zariski dense, then either $f$ is
homotopic to a holomorphic  map $f_0$, or there exists a compact
complex analytic space $Y'$, $\dim Y' < \dim Y$, and a holomorphic map
$\varphi  : Y\to Y'$, such that $f$ is homotopic to a composition
$Y\build \la^\varphi \fin Y' \build \la^{f_1}\fin  X$.

(4) Let $X = \Gamma \\ Sp(4)/U(2)$ be a Shimura threefold. Let $Y$
be compact K{\"a}hler and let $f : Y\to X$ be continuous. Then either $f$
is homotopic to a holomorphic map  $f_0$, or there exists a (singular) proper curve $S$, and a holomorphic map
$\varphi  : Y\to S$, such that $f$ is homotopic to a composition
$Y\build \la^\varphi \fin S \la X$.|

\rque Remarks|

1. I leave the case of Siegel type II to the reader (the proof is similar).

2. If $f_*$ is {\it not\/} rigid, one has strong consequences for $\pi
_1(Y)$, see 9.1. 

3. The lemma should be viewed as a final (twistorial) version of the
superrigidity theorem ([Si]).

\section { 2. Proof of the Superrigidity Lemma $(1)$}|

(1) Since we are given a continuous map $Y \build \la^f\fin X = \Gamma
\\ SU(m,n)/S(U(m)\times U(n))$ the map $S$ of Proposition~1.1 is
simply a holomorphic map $Y \to D/\Gamma $, where $D$ is a Griffiths
domain corresponding to the complex variation of Hodge structure,
defined by
$\rho  = f_* : \pi _1(Y) \to SU(m,n)$. Suppose the Higgs bundle looks
like $\oplus E^{p,q}$ where $E^{p,q}$ have dimensions
$m_1,n_1,m_2,n_2\ld m_s,k_s$ where $k_s$ is possibly missing. Then
$\summ m_i = m$, $\summ n_i=n$.
Now, the dimension of the horizontal distribution is
$$m_1\cdotp n_1 + n_1\cdotp m_2 + m_2 \cdotp n_2 + \cdots + m_s\cdotp
k_s~.$$
We notice that this number is strictly less than $m\cdotp n = \dim X$
except for the cases:

I)\hphantom{I}~~ $s=1$, \ie $E = E^{1,0} \oplus E^{0,1}$

II)~~ $s=2$, $k_2 =0$, \ie $E = E^{2,0} \oplus E^{1,1} \oplus
E^{0,2}$.

In the first case, $D$ is the symmetric space, and $D/\Gamma  = X$ so
we arrive to a holomorphic map to. $Y \to X$. In the second case the
second order equation reads $\theta _1(Z) \theta _2(W) - \theta
_1(W)\theta _2(Z) =0$ where $\theta _1 : TY \otimes E^{0,2} \to
E^{1,1}$ and $\theta _2 : TY \otimes E^{1,1} \to E^{2,0}$ are the
components of the (horizontal) derivative $DS$. So the image of $DS$
is strictly less than $\Hom(E^{0,2},E^{1,1}) \oplus \Hom
(E^{1,1},E^{1,0}) = m\cdotp n = \dim X$.  In other words, $\dim Y' <
\dim X$ where $Y' = S(Y)$.

\section { 3. Variation of Hodge structure, corresponding to rigid
representations to $SO(2,n)$}|

Let $\rho  : \pi _1(Y) \to SO(2,n)$ be a Zariski dense rigid
representation. Complexifying, we obtain a variation of Hodge structure
$E = \oplus E^{p,q}$. Since $\rho $ is defined over reals, we deal
with real variation of Hodge structure ([Sim~1]) that is to say,
$E^{p,q} = \overline{E^{p,q}}$ with respect to a flat complex
conjugation. For $n\ge 3$, this leaves exactly two possibilities:

I)\hphantom{I}~~ $E = E^{2,0} \oplus E^{1,1} \oplus E^{0,2}$, $\dim
E^{1,1} = n $, $\dim E^{2,0} = \dim E^{0,2} =1$.

II)~~ $n $ is even, $E = E^{2,0} \oplus E^{1,1} \oplus E^{0,2}$,
$\dim E^{1,1} =2$, $\dim E^{2,0} = \dim E^{0,2} = n/2$.

In case I) the Griffits domain is the symmetric spaces
$SO(2,n)/S(O(2) \times O(n))$ so the harmonic metric viewed as a
harmonic section of the flat bundle with fiber a symmetric space, is
holomorphic. In the second case the second order equation implies that
the rank of the derivative $DS$ of the $\rho $-equivariant holomorphic
map $\wt Y \to D$ is strictly less than $n$.

\section Proof of the Superrigidity Lemma {\rm (2)}|

This follows
immediately from the previous discussion in the same manner as in (1)

\section 4. Variations of Hodge structure, corresponding to rigid
representation to $Sp(4)$|

Let $\rho  : \pi _1(Y) \to Sp(2n)$ be a Zariski dense rigid
representation. Complexifying, we obtain a representation $\pi  : \pi
_1(Y) \to SU(n,n)$ and a real variation of Hodge structure $E =
\opp_{p+q=k} E^{p,q}$, $E^{p,q} = \overline E^{q,p}$ and $k$ odd. For
$n=2$ this leaves two possibilities:

\medskip
(I)~~ $E = E^{1,0} \oplus E^{0,1}$, and both $E^{1,0}$ and $E^{0,1}$,
or rather $V^{1,0}$ and $V^{0,1}$ viewed as $C^\infty$-subbundles of
the flat bundle $V$, are lagrangian wi
th respect to the flat complex
symplectic structure. This means first, that the Griffiths domain $D$
is the symmetric space $SU(2,2)/S(U(2)\times U(2)$,
second, that the image of the equivariant horizontal holomorphic map
$S : \wt Y \to D$ lies in the copy of the Siegel upper half-space
$Sp(4)/U(2)$ under the Satake embedding ([Sa]). In other words, the
unique $\rho $-equivariant harmonic map $\wt Y \to Sp(4)/U(2)$ is
holomorphic.

II)~~ $n $ is even, $E = E^{2,0} \oplus E^{1,1} \oplus E^{0,2}$,
$\dim E^{1,1} =2$, $\dim E^{2,0} = \dim E^{0,2} = n/2$.

In case I) the Griffits domain is the symmetric spaces
$SO(2,n)/S(O(2) \times O(n))$ so the harmonic metric viewed as a
harmonic section of the flat bundle with fiber a symmetric space, is
holomorphic. In the second case the second order equation implies that
the rank of the derivative $DS$ of the $\rho $-equivariant holomorphic
map $\wt Y \to D$ is strictly less than $n$.

\section Proof of the Superrigidity Lemma {\rm (2)}|

This follows
immediately from the previous discussion in the same manner as in (1)

\section 4. Variations of Hodge structure, corresponding to rigid
representation to $Sp(4)$|

Let $\rho  : \pi _1(Y) \to Sp(2n)$ be a Zariski dense rigid
representation. Complexifying, we obtain a representation $\pi  : \pi
_1(Y) \to SU(n,n)$ and a real variation of Hodge structure $E =
\opp_{p+q=k} E^{p,q}$, $E^{p,q} = \overline E^{q,p}$ and $k$ odd. For
$n=2$ this leaves two possibilities:

\medskip
(I)~~ $E = E^{1,0} \oplus E^{0,1}$, and both $E^{1,0}$ and $E^{0,1}$,
or rather $V^{1,0}$ and $V^{0,1}$ viewed as $C^\infty$-subbundles of
the flat bundle $V$, are lagrangian with respect to the flat complex
symplectic structure. This means first, that the Griffiths domain $D$
is the symmetric space $SU(2,2)/S(U(2)\times U(2)$,
second, that the image of the equivariant horizontal holomorphic map
$S : \wt Y \to D$ lies in the copy of the Siegel upper half-space
$Sp(4)/U(2)$ under the Satake embedding ([Sa]). In other words, the
unique $\rho $-equivariant harmonic map $\wt Y \to Sp(4)/U(2)$ is
holomorphic.

\medskip

(II)~~ $E = E^{3,0} \oplus E^{2,1} \oplus E^{1,2} \oplus E^{0,3}$ and
$\dim E^{p,q} =1$. The second order equation for $\theta $ implies
immediately that $D_s$ has rank at most one everywhere on $Y$.

\section Proof of the Superrigidity Lemma {\rm (4)}|

Follows from the discussion above.

\section { 6. Variations of Hodge structure, corresponding to rigid
representation to  $Sp(2n)$, and proof of the Superrigidity Lemma $(3)$}|

In general, the Higgs bundle is $E = \opp_{p+q=2s+1} E^{p,q}$,
$E^{p,q} = \overline {E^{p,q}}$. The dimension of the horizontal
distribution is 
$$ d = \summ_{p<s} \dim E^{p,q}\cdotp \dim E^{p+1,q-1}
= {\dim E^{s,s+1} \cdotp (\dim E^{s,s+1}+1)\over 2}~,$$
since $\theta  : E^{s,s+1} \to E^{s+1,s}$ viewed as bilinear form,
should be symmetric. Moreover, $\summ_{p\le s} \dim E^{p,q} = n$. An
elementary exercise shows that if $s>1$, $d < {n(n+1)\over 2}$. If
$s=1$, we get a holomorphic map to the Siegel upper half-plane.

\section 7. Regulators, I: proof of the Main Theorem|

The reader is supposed to be familiar with the geometric theory of
regulators ([Re~1], [Co]).

Let $\hb$ be a complex Hilbert space. The constant K{\"a}hler form
$(dX,dX)$ is invariant under the affine isometry group $\Iso(\hb)$, and $\hb$
is contractible, therefore there is a regulator class in
$H^2(\Iso^\delta (\hb),\rb)$. In fact, there is a class $\ell$ in
$H^1(\Iso^\delta (\hb),\hb)$ defined by a cochain $(x \mapsto Ux+b)
\longmapsto b$. The regulator class is simply $(\ell,\ell)$.

If $\pi _1(Y)$ does not have property $T$, then there exist a
representation $\rho :\pi _1(Y) \to \Iso(\hb)$ and a holomorphic
nonconstant section $S$ of the associated flat holomorphic affine
bundle with fiber $\hb$ ([Ko-Sch]). It follows that 
the pull-back $\rho ^*((\ell,\ell))$ of the regulator class to
$H^2(\pi _1(Y),\rb)$ restricts to a cohomology class in $H^2(Y,\rb)$,
given by a non-zero semi-positive $(1,1)$ form. Multiplying by the $\omega 
^{n-1}$, where $\omega $ is a K{\"a}hler form, and $n = \dim Y$, and
integrating over $Y$ we get a positive number, therefore this
cohomology class is non-zero. Therefore $H^2 (\pi _1(Y),\rb) \ne 0$.

Now if $\Delta$ does not have property T, and if $\pi _1(Y) \to \Delta$ is a central extension, then the construction of [Ko-Sch]  gives us an isometric uniform action on a real Hilbert space  $\pi _1(Y) \to \Iso(\hb)$, which factors through  $\Delta$, an
d
a harmonic section of an associated flat bundle. Since the action is uniform, the corresponding linear representation $\rho:\Delta \to  U(\hb)$ does not have fixed vectors. It follows from the  Lyndon-Serre-Hochschild spectral sequence that the map $H^1 (
\Delta ,\hb) \to H^1 (\pi _1(Y),\hb)$ is an isomorphism. By [Ko-Sch], thee exists an isometric action of $\pi _1(Y)$ on the complexified Hilbert space, extending the previous one, such that the corresponding flat bundle has a holomorphic section. This act
ion nesessarily factors through $\Delta$. Arguing as above, we deduce the theorem.

\rque Remark|Historically, the first break through in this direction
has been made in [JR], under assumption of having a nontrivial
variation of a finite-dimensional  unitary representation. Compare Proposition 9.1 below.

\rque Proof of the Corollary 0.2|The Lyndon-Serre-Hochschild spectral sequence implies that the map $H^2(
\Delta ,\rb) \to H^2(\wt
\Delta ,\rb )$ is zero. So $\wt \Delta $ is not a K{\"a}hler group.

\rque Remark|Suppose $\pi _1(Y)$ does not have property
$T$. Suppose moreover that that $\pi _1(Y)$ has a permutation
representation in $\ell^2(B)$, where $B$ is a countable set, and
$H^1(\pi _1(Y),\ell^2(B)) \ne 0$. Then we actually proved that
$H^2(\pi _1(Y),\ell^1(B)) \ne 0$. That is because the scalar product
$\ell^2(B) \times \ell^2(B) \to \cb$ factors through
$\ell^1(B)$. Moreover, the canonical map $H^2(\pi _1(Y), \ell^1(B))
\to H^2(\pi _1(Y),\cb)$ is nonzero.

\section 8. Regulators, II: proof of Theorems 0.1, 0.2|

Let $G$ be an isometry group of a classical symmetric bounded domain
$D$. With the exception of $SO(2,2)$, $H^1(G,\zb) = \zb$. This defines
a central extension $1\to \zb \to \wt G \to G \to 1$ and an extension
class  $e \in H^2(G^\delta ,\zb)$. On the other hand, the Bergman
metric on $D$ is $G$-invariant, so it defines a regulator class  $r
\in H^2_{\cont}(G,\rb)$.
It is proved in [Re~2], [Re~3] that, first, these classes coincide up
to a factor, and second, lie in the image of the bounded cohomology:
$H^2_b(G^\delta ,\rb) \to H^2(G^\delta ,\rb)$.

If $Y$ is a compact K{\"a}hler manifold, $\rho  : \pi _1(Y) \to G$ a representation , $s$ a holomorphic
nonconstant section of the associated flat $D$-bundle, then one sees
immediately that $(\rho ^*(r),\omega ^{n-1})>0$, so $\rho ^*(r)$, $\rho ^*(e)
\ne 0$. Theorem 0.1 follows now from the analysis of $VHS$ given in
sections 3, 4. To prove Theorem 0.2 notice that the case when $Y$ fibers over
a curve is obvious, otherwise $Y$ admits a holomorphic map to a
quotient of the Siegel half-plane and the proof proceeds as before.

\rque Remark|By [CT~1], the result of Theorem 0.1 is true for lattices
in $SU(n,1)$.

\section 9. Nonrigid representations|

\vskip-5pt
\th Proposition 9.1|Let $Y$ be compact K{\"a}hler and let $\rho  : \pi
_1(Y) \to SL(n,\cb)$ be a nonrigid irreducible representation. Then
$H^2(\pi _1(Y),\rb)\ne 0$.|

\preuve Let $\ggg = \slg(n,\cb)$ and let $\bar \rho $ be the adjoint
representation. We know that $H^1(\pi _1(Y),\ggg) \ne 0$. Therefore
$H^1(Y,\underline \ggg) \ne 0$ where $\underline \ggg$ is the local
system. By the Simpson's hard Lefschetz ([Sim~1]),
the multiplication by $\omega ^{n-1}$ gives an isomorphism
$H^1(Y,\underline \ggg) \to H^{2n-1} (Y,\underline \ggg)$, where
$\omega $ is the polarization class and $n = \dim_\cb Y$. The Poincar{\'e}
duality implies that the Goldman's pairing $H^1(Y,\ggg) \times
H^1(Y,\ggg) \to \cb$ is nondegenerate. Let $z$ be homology class in
$H_2(Y)$, dual to $\omega ^{n-1}$, and $\bar z$ its image in $H_2(\pi
_1(Y))$. It follows that the pairing $H^1(\pi _1(Y),\ggg) \times
H^1(\pi _1(Y),\ggg) \to \cb$ defined by $f,g\mapsto [(f,g),\bar z]$ is
nondegenerate. Here $(f,g) \in H^2(\pi _1,\cb)$ is the pairing defined
by the Cartan-Killing form. In particular, $\bar z \ne 0$.

\th Corollary 10.1|Let $Y$ a compact K{\"a}hler manifold. If $\pi _1(Y)$ has a Zariski dense
representation in either $Sp(4)$ or $SO(2,n)$, $n$ odd, then $H^2(\pi
_1(Y),\rb) \ne 0$.|

\preuve For rigid representations, this is proved in Theorems 0.1,
0.2. For nonrigid representations, this follows from Proposition~9.1.

\th Corollary 9.2|Let $\Gamma $ be any overgroup of a Zariski dense
countable subgroup of $Sp(4)$ or $SO(2,n)$, $n$ odd. Suppose
$b_1(\Gamma )=0$. Then the universal central extensions $\wt \Gamma $
is not K{\"a}hler.|

\section 10. Three-manifolds groups are not K{\"a}hler|

In this section, based on the previous development, we will present a
strong evidence in favour of the following conjecture, which we formulated in 1993 (Domingo Toledo informs us that a similar conjecture had been discussed by Goldman and Donaldson in 1989): 

\th Conjecture 10.1|Let $M^3$ be irreducible closed 3-manifold with
$\Gamma  = \pi _1(M)$ infinite. Then $\Gamma $ is not K{\"a}hler.|

\th Proposition 10.2 {\rm (Seifert fibration case)}|A cocompact
lattice in $\wt{SL(2,\rb)}$ is not K{\"a}hler.|

\preuve Passing to a subgroup of finite index, we can assume that
$\Gamma $ is a central extension of a surface group:
$$1 \la \zb \la \Gamma  \la \pi _1(S) \la 1$$
with a nontrivial extension class. In particular, $H^1(\Gamma ,\qb)
\simeq H^1(\pi _1(S),\qb))$, so the multiplication in $H^1(\Gamma
,\qb)$ is zero, which is impossible if $\Gamma $ is K{\"a}hler.

Recall that ``most'' of closed three-manifolds admit a Zariski dense
homomorphism $\pi _1(M) \build \la^\rho \fin SL_2(\cb)$ ([CGLS],
[Re~1]). 

\th Theorem 10.3|Let $M^3$ be atoroidal. Suppose there exists a
Zariski dense homomorphism $\rho  : \pi _1(M) \to SL_2 (\cb)$. Then
$\Gamma  = \pi _1(M)$ is not K{\"a}hler.|

\preuve By a theorem of [Zi]  $\pi _1(M)$ does not have
property $T$. By the Main Theorem, $H^2(\Gamma ,\rb) \ne 0$, hence by [Th],
$M$ is hyperbolic, which is impossible by [CT~1]. 

Alternatively, $\rho
$ is not rigid by [Sim~1], so $H^2(\Gamma ,\rb)\ne 0$ by Proposition
9.1, and then one procedes as before.

\rque Remark|In view of [CGLS], [Re~1], we obtain a huge number of
groups which are not K{\"a}hler.

\section 11. Central extensions of lattices in $PSU(2,1)$|

We saw  a general result, that, if $\Gamma  \subset SU(n,1)$ a
cocompact lattice and $[\omega ] \in H^2(\Gamma ,\zb)$ is given by any
ample line bundle, then a central extension
$$0 \la \zb \la \wt \Gamma  \la \Gamma  \la 1$$
with the extension class $[\omega ]$ is not K{\"a}hler. For $n=2$ one can
also prove:

\th Theorem 11.1|Let $\omega  \in H^2(B^2/\Gamma ,\zb) \cap (H^{2,0}\oplus
H^{0,2})$, $\omega \ne 0$. Then an extension
$$0\la \zb \la \wt \Gamma _\omega  \la \Gamma  \la 1$$
with the extension class $\omega $ is not K{\"a}hler.|

\rque Remark|$H^{2,0}$ becomes big on {\'e}tale finite coverings of
$B^2/\Gamma $ by Riemann-Roch.

\preuve Suppose $\wt \Gamma _\omega  = \pi _1(Y)$. The representation $\pi
_1(Y) \to \Gamma  \to SU(2,1)$ is rigid by the Lyndon-Serre-Hochschild
spectral sequence. It follows that there exists a dominating
holomorphic map $Y \to B^2/\Gamma $. But then the pullback map on
$H^{2,0}$ is injective, a contradiction.

\section 12. Smooth hypersurfaces in ball quotients which are not
$K(\pi ,1)$|

We saw that under various algebraic assumptions on $\Gamma  = \pi
_1(Y)$, there is a class in $H^2(Y,\rb)$ which vanishes on the
Hurewitz image $\pi _2(Y) \to H_2(Y,\zb)$, therefore defining a
nontrivial element of $H^2(\Gamma ,\rb)$. On the contrary, we will
show now  that there are hypersurfaces $Y$ in ball quotients
$B^n/\Gamma $, $n\ge 3$ with a surjective map $\pi _1(Y) \to \Gamma $
such that $\pi _i(Y) \ne 0$ for some $i$. The proof is very indirect
and we don't know the exact value of $i$. The varieties $Y$ were in fact
introduced in [To] where it is proved that $\pi _1(Y)$ is not
residually finite. We will show that $cd(\pi _1(Y)) \ge 2n-1$,
therefore $Y$ is not $K(\pi ,1)$.

Let $X^n$ be an arithmetic ball quotient and let $X_0 \subset X$ be a
totally geodesic smooth hypersurface. Let $D = X - X_0$,
 then $D$ is covered topologically by $\cb^n$ minus a countable
union of hyperplanes, so $D$ is $K(\pi ,1)$. Let $S$ be a boundary of
a regular heighbourhood of $X_0$, so $S$ is a circle bundle over
$X_0$, in particular $S$ is $K(\pi ,1)$ and $\pi _1(S)$ is a central
extension $0 \to \zb \to \pi _1(S) \to \pi _1(X_0) \to 1$ with a
nontrivial extension class (this is because the normal bundle to $X_0$
is negative). Let $V$ be a  finite dimensional module over $\pi _1(X)$
with an invariant nondegenerate form $V\to V'$. We have an exact sequence
$$\displaylines{
H^{2n-1}(\pi _1(X),V) \la H^{2n-1}(\pi _1(X_0),V) \oplus H^{2n-1}(\pi
_1(D),V) \hfill\cr
\hfill \la H^{2n-1}(\pi _1(S),V)
 \la H^{2n}(\pi _1(X),V) \la \cdots}$$
Now, we make a first assumption:

\medskip
1) $H^0(\pi _1(X),V)=0$. 

It follows that $H_0 (\pi _1(X),V) =0$, so
$H^{2n}(\pi _1(X),V) =0$; we make a second assumption : 
\medskip

2) $H^1(\pi _1(X),V) =0$. 

It follows that
$H^{2n-1}(\pi _1(X),V) =0$. So we have (remember that $X_0$ has
dimension $n-1$)
$$H^{2n-1}(\pi _1(D),V) \simeq H^{2n-1} (\pi _1(S),V)~.$$
Now, in the $E^2$ of the Lyndon-Serre-Hochschild spectral sequence for
$H^*(\pi _1(S),V)$ the term $H^{2n-2}(\pi _1(X_0),H^1(\zb,V))$ is not
hit by any differential. Since $\zb$ acts trivially, this is just
$H^{2n-2}(\pi _1(X_0),V) \simeq H_0(\pi _1(X_0),V)$. We now make a
third assumption:

\medskip
3) $H^0(\pi _1(X_0),V) \ne 0$.

Then we will have $H^{2n-1}(\pi _1(D),V) \ne 0$. Let $Y$ be a generic
hyperplane section of $X/X_0$, constructed in [To], then [GM], $\pi
_1(Y) = \pi _1(D)$ and we are done.

Now, we take for $V$ the adjoint module. The assumption 2) follows
from Weil's rigidity. The assumption 3) is satisfied for standard
examples of $X_0$ ([To]).

\rque Remark|The construction of [To] is given for lattices in
$SO(2,n)$, but it applies verbatim here.

\vfill\eject
\vphantom{}

\vskip-2cm

\titre References|

{\ept
\lettre Gr-Sch|

\article Co|K. Corlette|Rigid representations of K{\"a}hlerian fundamental
groups|J.D.G.|33|1991|239--252|

\article CT~1|J. Carlson, D. Toledo|Harmonic maps of K{\"a}hler manifolds to
locally symmetric spaces|Publ. IHES|69|1989|173--201|

\article CGLS|M. Culler C. Gordon, J. Luecke, P. Shalen|Dehn surgery on
knots|Annals of Math.|{\bf 125}|1987|237--308|


\livre G|Ph. Griffiths, ed.|Topics in transcendental algebraic
geometry|Princeton UP|1984|

\article Gr-Sch|M. Gromov, R. Schoen|Harmonic maps into singular
spaces and $p$-adic superrigidity for lattices in groups of rank
one|Publ. Math. Inst. Hautes {\'E}tud. Sci.|76|1992|165--246|

\article  GS|Ph. Griffiths, W. Schmidt|Locally homogeneous complex
manifolds|Acta Math.|123|1969|253--302|

\livre   H|P. {\rm de la} Harpe, A.Valette|La propri\'et\'e (T) de kazhdan pour les groupes localement compactes|Ast\'erisque|1989|

\article Hit|N. Hitchin|The self-duality equation on a Riemann
surfaces|Proc. London Math. Soc.|55 (3)|1987|59--126|

\article  JR|F.E.A.Johnson, E.G.Rees|On the fundamental group of a complex algebraic manifold|Bull. Lond. Math. Soc.|19|1987|463--466|

\livre  Ko|J. K\`ollar|Shafarevich maps and automorphic forms|Princeton UP|1995|

\article  Ko-Sch|N. Korevaar, R. Schoen|Global existence theorems for
harmonic maps to non locally compact spaces|Commun. Anal. Geom.|5
(2)|1997|333--387|

\article Re~1|A. Reznikov|Rationality of secondary
classes|J. Differential Geom.|43|1996|674--692|

\article Re~2|A. Reznikov|Three-manifolds class field theory|Selecta
Math.|3|1997|361--399|

\livre  Sa|I. Satake|Algebraic structures of symmetric
domains|Princeton UP|1980|

\article Sim~1|C. Simpson|Higgs bundles and local
systems|Publ. IHES|75|1992|5--95|

\article Sim~2|C. Simpson|Moduli of representations of the fundamental
group of smooth projective variety|Publ. IHES|79-80|1994|47--129, 5--79|

\article Si|Y.-T. Siu|The complex-analyticity of harmonic maps and the
strong rigidity of compact  K{\"a}hler manifolds|Ann. Math.|112|1980|73--111|

\livre Th|W. Thurston|The geometry and topology of three
manifolds|Lecture notes|{\rm Princeton University}|

\article To|D. Toledo|Projective varieties with non-residually finite
fundamental group|Publ. IHES|77|1993|103--119|

\article  Zi|R. Zimmer|Kazhdan groups acting on compact manifolds|Invent. Math.|75|1984|425--436|
}
\vskip-0.5cm
$$-\diamondsuit-$$
\vfill
{\obeylines\parskip=0pt \ept
Department of Mathematical Sciences
University of Durham
South Road
Durham DH1 3LE(UK)
\medskip
{\it Current address~:}
Universit{\'e} de Grenoble I
{\bf Institut Fourier}
UMR 5582
UFR de Math{\'e}matiques
B.P. 74
38402 St MARTIN D'H{\`E}RES Cedex (France)
\medskip
Alexander.Reznikov@durham.ac.uk
}

\end